\documentclass{amsart}

\usepackage{amssymb}
\usepackage{url}

\newtheorem{theorem}{Theorem}[section]
\newtheorem{lemma}[theorem]{Lemma}
\newtheorem{proposition}[theorem]{Proposition}

\newcommand\spc{\mkern+1mu}
\newcommand\nspc{\mkern-1mu}

\numberwithin{equation}{section}

\def\lra{\longrightarrow}
\DeclareMathOperator{\Aut}{Aut}

\begin{document}

\title[On the nonexistence of Smith-Toda complexes]
{On the nonexistence of Smith-Toda complexes}

\author{Lee S. Nave}
\address{University of Washington, Department of Mathematics,
Box 354350, Seattle, WA  98195--4350}
\email{nave@math.washington.edu}

\date{October 30, 1998}

\begin{abstract}
Let $p$ be a prime.  The Smith-Toda complex $V(k)$ is a finite spectrum
whose $BP$-homology is isomorphic to $BP_{*}/(p,v_{1},\dots,v_{k})$.  For
example, $V(-1)$ is the sphere spectrum and $V(0)$ the mod~$p$ Moore spectrum.
In this paper we show that if $p>5$, then $V((p+3)/2)$ does not exist and
$V((p+1)/2)$, if it exists, is not a ring spectrum.  The proof uses the new
homotopy fixed point spectral sequences of Hopkins and Miller.
\end{abstract}

\maketitle

\section{Introduction}

Let $p$ be a prime and recall the Brown-Peterson spectrum $BP$.  It is
a $p$-local ring spectrum with coefficient ring $BP_*\cong\mathbb
Z_{(p)}[v_1,v_2,\dots]$, where $v_i$ is the $i$th Hazewinkel generator
in degree $2(p^i-1)$.  If $X$ is a spectrum, $BP_*X$ is a comodule
over the Hopf algebroid $BP_*BP$.  (See \cite{rav:green} for details.)
Smith~\cite{smi:onreal} and 
Toda~\cite{tod:onspec}
considered the existence of finite spectra $V(k)$ with $BP_*V(k)\cong
BP_*/I_{k+1}$ (as $BP_*$-modules, hence as $BP_*BP$-comodules), where
$I_{k+1}=(p,v_1,\dots,v_k)$.  For example, $V(-1)=S^0$ and
$V(0)=M(p)$, the mod~$p$ Moore spectrum.

The construction of $V(1)$ is originally due to Adams~\cite{ada:jx4}.
Smith constructed $V(2)$ for $p>3$ and Toda constructed
$V(3)$ for $p>5$.  Furthermore, these results are sharp.
The first negative results were obtained by Toda, who showed that
$V(1)$ cannot be constructed when $p=2$ and likewise for $V(2)$ when
$p=3$.  Later, Ravenel~\cite[7.5.1]{rav:green} showed that $V(3)$
does not exist when $p=5$.

Recently, Hopkins, Mahowald, and Miller established the nonexistence
of $V(p-2)$ for all $p>3$.  In this paper, we prove

\begin{theorem}\label{thm:main}
If $p>5$, $V((p+3)/2)$ does not exist.  If $V((p+1)/2)$ exists, it
is not a ring spectrum.
\end{theorem}

The proof uses consequences of new work of Hopkins and Miller, not yet
published, which we will briefly describe.
An account of this work
is given in \cite{rez:noteso}.
Throughout this paper, we restrict ourselves to the 
category of $p$-local spectra.  We call a spectrum $p$-compact if it is
the $p$-localization of a finite spectrum.  (The $p$-compact spectra
comprise the small objects in the $p$-local stable homotopy category,
in the sense of~\cite{hov-pal-str:axioma}.)  We write
$S^0$ for the $p$-local sphere spectrum.

First we must recall the relevant parts of Morava's
theory (\cite{mor:noethe}, see also \cite{dev:morava}).  Fix
$n\geq1$, set $q=p^{n}$,
and let $E_{n*}\equiv W\mathbb F_{q}[\nspc[u_{1},\dots,
u_{n-1}]\nspc][u,u^{-1}]$, with grading
$|u|=-2$ and $|u_i|=0$.  Here $W\mathbb F_{q}$ denotes
the Witt vectors with coefficients in $\mathbb F_{q}$.  $E_{n*}$ is
a graded $BP_*$-algebra via the map $\theta:BP_*\to E_{n*}$ given by
\begin{equation*}
\theta(v_i) = \begin{cases} u^{1-p^i}u_i,&i<n\\
                            u^{1-p^n},   &i=n\\
                            0,           &i>n.
              \end{cases}
\end{equation*}
We will write $v_i$ for $\theta(v_i)\in E_{n*}$.
By the Landweber exact functor theorem \cite{lan:homolo},
$E_{n*}\otimes_{BP_*}(-)$ is a homology theory and we write $E_n$
for the representing spectrum.  If $X$ is a $p$-compact spectrum, Morava's
theory provides an action of $S_n$ on $E_{n*}X$, where $S_n$ is
the $n$th Morava stabilizer group.  Furthermore, 
$\text{Gal}\equiv\text{Gal}(\mathbb F_{q}/\mathbb F_p)$ acts on
$E_{n*}X$ via its action on $W\mathbb F_{q}$, making 
$E_{n*}X$ into a $G_n$-module, where $G_n\equiv S_n\rtimes\text{Gal}$.
In fact, this action is induced by an action (in the homotopy category)
of $G_n$ on $E_n$ (\cite{dev:moravam}, see discussion p.\ 767).

Morava's change of rings theorem, in conjunction with
\cite{mah-sad:vntele}, leads to a spectral sequence
\begin{equation*}
H^*_c(G_n;E_{n*}X)\Longrightarrow \pi_*L_{K(n)}X,
\end{equation*}
where $L_{K(n)}$ denotes localization with respect to $K(n)$,
the $n$th Morava $K$-theory spectrum.
This spectral sequence resembles a homotopy fixed point spectral
sequence, but with two crucial differences.  First, the group $G_n$ is
not discrete.  In fact, it is profinite and the $E_2$ term involves
continuous cohomology.  Second, $G_n$ acts on $E_n$ only up to homotopy.
To form homotopy fixed point spectra in the usual way
requires that the action be ``on the nose''.

Nevertheless, Hopkins and Miller have shown that if
$G\subset G_n$ is a finite subgroup one may form the homotopy fixed
point spectrum $E_n^{hG}$ and if $X$ is a $p$-compact spectrum there is a
spectral sequence
\begin{equation*}
H^*(G;E_{n*}X)\Longrightarrow \pi_* (E_n^{hG}\wedge X).
\end{equation*}
Subsequent work of Devinatz and Hopkins, in preparation, has shown
that this can be done for any closed subgroup $G\subset G_n$, provided
one works with continuous cohomology.
More precisely, the spectral sequence is the $K(n)$-local
$E_n^{\text{Gal}}$-based Adams spectral sequence for 
$E_n^{hG}\wedge X$, where
$E_n^{\text{Gal}}$ is the Landweber exact spectrum with
coefficient ring $\mathbb Z_p[\nspc[ u_1,\dots,u_{n-1}]\nspc][u,u^{-1}]$.

This construction is natural in $G$ and agrees with the
usual spectral sequence when $G=S_n$.  Also, the spectral sequence is
multiplicative when $X$ is a ring spectrum.  A suitable geometric
boundary theorem holds as well:
Suppose
\begin{equation*}
X\overset{f}{\lra} Y\overset{g}{\lra} Z\overset{h}{\lra}\Sigma X
\end{equation*}
is a cofiber sequence of $p$-compact spectra with
$E_{n*}(h)=0$.  The resulting short exact sequence
\begin{equation*}
0\lra E_{n*}X\lra E_{n*}Y\lra E_{n*}Z\lra0
\end{equation*}
gives rise to a connecting homomorphism
$H^*_c(G;E_{n*}Z)\overset{\delta}{\to} H^{*+1}_c(G;E_{n*}X)$.
Suppose $z\in\pi_*(E_n^{hG}\wedge Z)$ is
detected by $\bar z\in H^*_c(G;E_{n*}Z)$.  Then $\delta(\bar z)$
is a permanent cycle and
$h(z)\in
\pi_{*-1}(E_n^{hG}\wedge X)$ is detected by
$\delta(z)$, or else $h(z)$ is in higher filtration.

For the rest of the paper, we take $p$ odd and $n=p-1$.  In this case,
$S_n$ has exactly one maximal finite subgroup $G$, unique up to
conjugation, namely $G\cong\mathbb Z/p\rtimes\mathbb Z/n^2$, where the
action of the right factor on the left is given by
\begin{equation*}
\mathbb Z/n^2\to\mathbb Z/n\cong\Aut(\mathbb Z/p).
\end{equation*}
(See \cite{hew:finite}, \cite{hew:normal}, also \cite{mil:noteso}.)
If $X$ is a $p$-compact spectrum, we refer to the spectral sequence
\begin{equation*}
H^*(G;E_{n*}X)\Longrightarrow \pi_*(E_n^{hG}\wedge X)
\end{equation*}
as the {\bf spectral sequence for $X$}.

\section{Preliminaries}

We begin with a discussion of the spectral sequence for $S^0$.  These
results are due to Hopkins and Miller.  Because they have
yet to appear in print, we include some details for the reader's
convenience.  Full proofs will also appear in the author's thesis.  We
remind the reader that $n=p-1$, $q=p^{n}$, and $G\cong\mathbb Z/p\rtimes
\mathbb Z/n^2$.

Our interest in this work began when Hopkins showed us how to
compute the $E_2$ term of this spectral sequence.  To wit,
\begin{equation}\label{eq:cohomology}
H^*(G;E_{n*})\cong \mathbb F_q[\Delta^{\pm}]\langle\alpha
\rangle[\spc\beta\spc],
\end{equation}
where $|\Delta|=(0,2pn^2)$, $|\alpha|=(1,2n)$, and $|\beta|=(2,2pn)$.
(The cohomological degree is given first.)  Actually, this
presentation ignores much of $H^0(G;E_{n*})$.  More precisely,
$H^0(G;E_{n*})/p$ has a splitting with $\mathbb F_q[\Delta^\pm]$
as a direct summand.  With this splitting, the projection
$H^0(G;E_{n*})\to\mathbb F_q[\Delta^\pm]$ gives $\mathbb
F_q[\Delta^\pm]\langle\alpha\rangle[\spc\beta\spc]$ the structure of an
$H^0(G;E_{n*})$-algebra, which is isomorphic to $H^*(G;E_{n*})$ in
positive cohomological degrees.

The calculation of $H^*(G;E_{n*})$ involves choosing coordinates for
$E_{n*}$ on which the action of $G$ is easily described.  Let
$\mathfrak m =(p,u_1,\dots,u_{n-1})$ be the maximal ideal of $E_{n*}$.
There are elements $w$, $w_1$,\dots,~$w_{n-1}\in E_{n*}$ with
\begin{alignat*}{2}
w&\equiv u&&\mod (p,\mathfrak m^2)\\
  w_i&\equiv u_i&&\mod (p,u_1,\dots,u_{i-1},\mathfrak m^2)
\end{alignat*}
and $G$ action given by
\begin{align*}
\sigma(w) &= w+ww_{n-1},\\
\sigma(ww_i) &= ww_i + ww_{i-1},\quad\text{for $2\leq i\leq n-1$,}\\
\tau(w)&=\eta w,
\end{align*}
and the relation $(1+\sigma+\dots+\sigma^{p-1})w=0$,
where $\sigma$ generates $\mathbb Z/p$, $\tau$ generates $\mathbb
Z/n^2$, and $\eta$ is a primitive $n^2$ root of unity.

To obtain this representation of $E_{n*}$ requires the following
formulae for the action of $G$ on the generators $u$, $u_i$, which
we also record for later use.  Recall that an element $g\in S_n$
has a unique expression of the form $g=\sum_{j=0}^{\infty} a_j S^j$,
where $a_j\in W\mathbb F_{q}$, $a_j^{q}=a_j$, and $a_0$ is a unit.
We have
\cite[Proposition~3.3 and Theorem~4.4]{dev-hop:theact}
\begin{equation*}
  g(u)\equiv a_0u+a_{n-1}^{\chi}uu_1+\dots+a_1^{\chi^{n-1}}
  uu_{n-1}\mod(p,\mathfrak m^2) 
\end{equation*}
and
\begin{equation*}
  g(uu_i)\equiv a_0^{\chi^i}uu_i+\dots+a_{i-1}^{\chi}uu_1
  \mod(p,\mathfrak m^2),
\end{equation*} 
where $\chi$ denotes the
Frobenius automorphism of $W\mathbb F_{q}$.  Now it turns out that
$\sigma=\sum_{j=0}^{\infty}a_j S^j$ with $a_0\equiv1$ mod~$p$ and
$a_1\in (W\mathbb F_{q})^{\times}$.  Also, $\tau=\eta$.  Then for
$1\leq i\leq n$,
\begin{equation}\label{eq:action}
  \begin{aligned}
    \sigma(uu_i) &\equiv uu_i+c_iuu_{i-1}&&\mod
     (p,u_1,\dots,u_{i-2},\mathfrak m^2)\\
    \tau(uu_i) &\equiv\eta^{p^i}uu_i&&\mod(p,\mathfrak m^2),
  \end{aligned}
\end{equation}
where $c_i\in (W\mathbb F_{q})^\times$ and we make the conventions
that $u_0=p$ and $u_n=1$.

To get the desired representation of $E_{n*}$, Hopkins and Miller
start with
$s=(1-\sigma)(v_1)/p\in E_{n*}$.
From the definition of the action of $S_n$ on $E_{n*}$ 
\cite[(5.2)]{dev:morava} and the formula $\eta_R(v_1)=v_1+pt_1$ in 
$BP_*BP$, we have $s=t_1(\sigma^{-1})$, where we write $t_i$ for 
the image of $t_i$ under the map of Hopf algebroids
\begin{equation*}
(BP_*,BP_*BP)\lra(E_{n*}^{\text{Gal}},\text{Map}_c(S_n,E_{n*})^{\text{Gal}}).
\end{equation*}
This map is a composition of several Hopf algebroid maps.
By unraveling
the definitions of these maps in \cite{dev:morava}, it is clear that
if $g\in S_n$ is given by $g=\sum_{j=0}^{\infty}a_j S^j$, then
$t_i(g)\equiv a_0^{-1}a_i u^{1-p^i}$ mod~$\mathfrak m$.

In particular, $s\equiv c u^{1-p}$ mod~$\mathfrak m$, where 
$c\in(W\mathbb F_{q})^{\times}$.  Let 
\begin{equation*}
t=s\prod_{j=0}^{p-1}\sigma^j(u)
\qquad\text{and}\qquad
w=\frac{1}{n^{2}}\sum_{j=1}^{n^2}\eta^{-j}\tau^j(t).
\end{equation*}
Then 
$w\equiv cu$ mod~$(p,\mathfrak m^2)$,
$(1+\sigma+\dots+\sigma^{p-1})w=0$, and
$\tau(w)=\eta w$.
Finally, let $w_i=(\sigma-1)^{n-i}(w)$.

The differentials are determined by the Toda differential in the
Adams-Novikov spectral sequence (ANSS) and the nilpotence of
$\beta_1\in\pi_*S^0$ as follows.
There is a map from the ANSS to the spectral sequence
for $S^0$,
\begin{equation*}
\text{Ext}_{BP_*BP}(BP_*,BP_*)\lra H^*(G;E_{n*}).
\end{equation*}
By \cite{rav:thenon}, this map sends
$\alpha_1$ to $\alpha$, $\beta_1$ to $\beta$, and
$\beta_{p/p}$ to $\Delta\beta$.  (We will omit the phrase
``up to multiplication by a unit'' from this discussion.)
(See \cite[Section 1.3]{rav:green} for notation.)
Because $\alpha_1$ and $\beta_1$ are permanent cycles, so are $\alpha$
and $\beta$.

Let $K$ denote the kernel of the projection $H^{0}(G;E_{n*})\to
\mathbb F_{q}[\Delta^{\pm}]$.  Because $\beta\cdot K=0$, $d_{2p-1}$
vanishes on $K$.
Also, from the Toda differential
$d_{2p-1}(\beta_{p/p})=\alpha_1\beta_1^p$
(\cite{tod:animpo},\cite{tod:extend}) and
the evident sparseness in our spectral sequence, we conclude that
\begin{equation*}
d_{2p-1}(\Delta)=\alpha\beta^{p-1}.
\end{equation*}
These facts and the multiplicative structure of the
spectral sequence completely determine $d_{2p-1}$.

For degree reasons the next
possible differential is $d_{2n^2+1}$.  Since $\beta_1^{pn+1}=0$ in
$\pi_*S^0$ (\cite{tod:animpo}), \cite{tod:extend}), $\beta^{pn+1}$ is
hit by some differential in our spectral sequence.  The last
differential which could do this is $d_{2n^2+1}$ and therefore
\begin{equation*}
d_{2n^2+1}(\Delta^{p-1}\alpha)=\beta^{n^2+1}.
\end{equation*}
Again, $d_{2n^{2}+1}$ vanishes on $K$.  This determines $d_{2n^2+1}$
completely and the spectral sequence collapses after this point for
degree reasons, hence

\begin{proposition}\label{prop:infinity}
If $s$ is odd, then $E_\infty^{s,t}=0$ unless $1\leq s\leq 2n-1$ and
\begin{equation*}
t\equiv 2n+(s-1)pn+2pn^2x\quad\text{mod $2p^2n^2$},\quad
\text{where $x\not\equiv-1$ mod $p$.}
\end{equation*}
If $s>0$ is even, then $E_{\infty}^{s,t}=0$ unless $2\leq s\leq 2n^2$ and
\begin{equation*}
t\equiv spn\quad\text{mod $2p^2n^2$.}
\end{equation*}
\end{proposition}

We also need a handle on the zero line.
Using the coordinates $w$, $w_i$ mentioned above, it is easy to show
that

\begin{proposition}\label{prop:fixed}
$H^0(G;E_{n*})$ is concentrated in degrees $t\equiv0$ mod~$2n$.
\end{proposition}

Finally, we will need the following partial description of
$H^*(G;E_{n*}V(1))$, which is a simple consequence of~\eqref{eq:cohomology}.

\begin{proposition}\label{prop:v1}
Let $A$ be the graded abelian group given by
\begin{equation*}
A_m=H^{1,2pn+2pn^2m}(G;E_{n*}/I_2),
\end{equation*}
made into
an $\mathbb F_{q}[\Delta^\pm]$-module via the map
$E_{n*}\to E_{n*}/I_2$.  Then
$A\cong\Sigma^{2pn}\mathbb F_{q}[\Delta^\pm]$.
\end{proposition}

\section{Proof of Theorem~\ref{thm:main}}

Suppose $V(k)$ exists.  It is well known that the condition
$BP_*V(k)\cong BP_*/I_{k+1}$ is equivalent to 
$H_*(V(k);\mathbb F_p)$ being isomorphic to the exterior algebra
$E(\tau_0,\dots,\tau_k)$ as a comodule over the mod~$p$ dual
Steenrod algebra. 
This point of view allows one to produce
$V(k-1)$ as a skeleton of $V(k)$ and exhibit $V(k)$ as the cofiber
of a map $f:\Sigma^{|v_k|}V(k-1)'\to V(k-1)$ inducing
multiplication by $v_k$ in $BP$-homology.  
(The notation $V(k-1)'$ is meant to distinguish this spectrum from
$V(k-1)$---they need not be isomorphic.)

Recall that in the 
Hopf algebroid structure of $BP_*BP$, $v_k$ is invariant mod~$I_k$,
i.e., $\eta_r(v_k)\equiv v_k$ mod~$I_k$.  Then 
$v_k\in E_{n*}V(k-1)\cong E_{n*}/I_k$ is fixed by all of $S_n$.
In particular, $v_k^j\in H^0(G;E_{n*}V(k-1))$ for all $j\geq0$.

The map $f$ induces a map (of degree $|v_k|$) from the spectral sequence
for $V(k-1)'$ to the spectral sequence for $V(k-1)$, which is multiplication
by $v_k$ on $H^0$.  In particular, if $v_k^j$ is a permanent cycle 
in the spectral sequence for $V(k-1)'$, then $v_k^{j+1}$ is a
permanent cycle in the spectral sequence for $V(k-1)$.

Fix $p>5$ and let $m=(p+3)/2$.  We show that if $V(m-1)$ exists,
then in the spectral sequence for $V(m-1)$, $v_m$ is a permanent cycle
and $v_m^2$ is not.  Then clearly $V(m-1)$ cannot be a ring spectrum,
and by the discussion above, $V(m)$ cannot exist.  In order to
simplify the notation, we will ignore the distinction between
$V(k-1)'$ and $V(k-1)$.  Since our arguments depend only on properties
of the spectral sequence for $S^0$, which is unique in the $p$-local
category, no generality is lost.

Before proceeding, we need to establish the vanishing of
$H^{s,t}(G;E_{n*}/I_{k})$ and $\pi_r(E_n^{hG}\wedge V(k))$
for various values of $s$, $t$, $r$, and~$k$.  We are not able
to calculate these groups directly for arbitrary values of $k$.
Instead, we use a brute force approach which reduces the
computation to the case $k=-1$.

\begin{lemma}\label{lem:cohomology} 
If $d\in\mathbb Z$ satisfies
\begin{equation*}
d\equiv 2pny\mod 2pn^2
\end{equation*}
where $y\not\equiv 1$ mod~$n$, then
\begin{equation*}
H^{1,d}(G;E_{n*}/I_k)=0.
\end{equation*}
If, in addition, $k\ne p-1$, then
\begin{equation*}
H^{2,d-2n}(G;E_{n*}/I_k)=0.
\end{equation*}
\end{lemma}

\begin{proof}
If $k\geq p$, $E_{n*}/I_k=0$ so the result is trivial.  Thus assume
$k<p$.
We prove first
that $H^{1,d}(G;E_{n*}/I_k)=0$.  Using
the long exact sequences in cohomology arising from the short exact
sequences of $E_{n*}$-modules
\begin{equation*}
0\lra E_{n*}/I_i\overset{v_i}{\lra} E_{n*}/I_i\lra E_{n*}/I_{i+1}\lra 0,
\end{equation*}
it suffices to show that $H^{s,t}(G;E_{n*})=0$ for all $(s,t)$ with
$1\leq s\leq k+1$ and
\begin{equation}
t=d-\sum_{j=1}^{s-1}|v_{i_j}|,\label{eq:sum1}
\end{equation}
where $0\leq i_1<i_2<\dots<i_{s-1}\leq k-1$.

First, suppose $s$ is even.  By \eqref{eq:cohomology},
$H^{s,t}(G;E_{n*})=0$ unless
$t\equiv spn$ mod~$2pn^2$.  Substituting this into~\eqref{eq:sum1} and
reducing mod~$p$ leads to $s=2$ and $i_1=0$.  In this case,
working mod~$2pn^2$ leads to $y\equiv1$ mod~$n$, contrary to our
hypothesis on $y$.

Now suppose $s$ is odd.  Again by \eqref{eq:cohomology},
$H^{s,t}(G;E_{n*})=0$ unless
$t\equiv 2n+(s-1)pn$ mod~$2pn^2$.  Substituting this into~\eqref{eq:sum1}
and reducing mod~$p$ leads to an equation in $s$ which, keeping in mind
all the constraints involved, has no solution.

Similarly, we reduce the proof that
$H^{2,d-2n}(G;E_{n*}/I_k)=0$ to showing that $H^{s,t}(G;E_{n*})=0$
for all $(s,t)$ with $2\leq s\leq k+2$ and
\begin{equation}
t=d-2n-\sum_{j=1}^{s-2}|v_{i_j}|,\label{eq:sum2}
\end{equation}
where $0\leq i_1<i_2<\dots<i_{s-2}\leq k-1$.

As before, we proceed by reducing~\eqref{eq:sum2} mod~$p$.  This 
time there are no solutions for $s$, assuming as we are that $k<p-1$.
\end{proof}

\begin{lemma}\label{lem:homotopy}
Let $m<p-1$.  If $d\in\mathbb Z$ satisfies
\begin{equation*}
d\equiv 2n+2pn+2p^2ny\mod 2p^2n^2,
\end{equation*} 
where $y\not\equiv0$ mod~$n$, and $V(m)$ exists, then
\begin{equation*}
\pi_{d-1}(E_n^{hG}\wedge V(m))=0.
\end{equation*}
\end{lemma}

\begin{proof}
Using the long exact sequences in homotopy arising from the cofiber
sequences
\begin{equation*}
\Sigma^{|v_i|}V(i-1)\overset{v_i}{\lra} V(i-1)\longrightarrow V(i),
\end{equation*}
it suffices to show that in the spectral sequence for $S^0$,
$E_{\infty}^{s,t}=0$ for all $(s,t)$ with
\begin{equation}
t-s=d-1-k-\sum_{j=1}^k |v_{i_j}|,\label{eq:sum3}
\end{equation}
where $0\leq k\leq m+1$ and $0\leq i_1<i_2<\dots<i_k\leq m$.

First, suppose $s=0$.  Reducing~\eqref{eq:sum3} mod~$2n$, we have $t\equiv
-(k+1)$.  By Proposition~\ref{prop:fixed}, $E_{\infty}^{0,t}=0$.

Next, suppose $s$ is odd.  By Proposition~\ref{prop:infinity}, we may
take $1\leq s\leq 2n-1$ and $t\equiv 2n+(s-1)pn$ mod~$2pn^2$.
Substituting this into~\eqref{eq:sum3} and reducing mod~$2n$, we are
led to $s=k+1$.  Substituting this into~\eqref{eq:sum3} and reducing mod~$p$
yields $k=0$ or $k=1$ (if $i_1=0$).  The second case cannot occur because
$s$ is odd, therefore we may take $k=0$ and $s=1$.  In this case,
\eqref{eq:sum3} becomes $t=d$.  Reducing this mod~$2pn^2$, we are led to
$y\equiv-1$ mod~$n$.
It follows that
$t\equiv 2n-2pn^2$ mod~$2p^2n^2$.  By Proposition~\ref{prop:infinity},
$E_{\infty}^{1,t}=0$.

Finally, suppose $s>0$ is even.  In this case, we may take $2\leq s\leq 2n^2$
and $t\equiv spn$ mod~$2p^2n^2$.  Working mod~$2n$ as above, we get
$s=2nl+k+1$, where $0\leq l<n$.  Note that $k$ must be odd.
Substituting these equations for $s$ and $t$ into~\eqref{eq:sum3} and 
reducing mod~$p$ leads to $l=k-1$ (if $i_1>0$) or $l=k-2$ (if $i_1=0$).
To finish the proof, we need to reduce~\eqref{eq:sum3} mod~$2p^2n$, which is
a bit unwieldy.  There are four cases:
\begin{equation*}
t-s\equiv 2n+2pn-(k+1)-\cases (2pn+2n)k,&\text{$i_1>1$}\\
	 (2pn+2n)(k-1)+2n,&\text{$i_1=1$}\\
	 (2pn+2n)(k-1),&\text{$i_1=0$ and $i_2>1$}\\
	 (2pn+2n)(k-2)+2n,&\text{$i_1=0$ and $i_2=1$.}
\endcases
\end{equation*}

For each case, we substitute $t=spn$, $s=2nl+k+1$, and either $l=k-1$
(first two cases) or $l=k-2$ (last two cases) and solve for $k$.  The
first and third cases yield $k\equiv-1$ mod~$2p$, which is not
possible.  The second and fourth cases yield $k\equiv1$ mod~$2p$,
i.e., $k=1$.  In the fourth case, this is not possible, because we are
assuming that $i_2=1$ in that case, so $k>1$.  This leaves the second
case, with $k=1$.  Substituting $t=spn$, $s=2nl+k+1$, $l=k-1$, and
$k=1$ into~\eqref{eq:sum3} and reducing mod~$2p^2n^2$ leads to
$y\equiv0$ mod~$n$, contrary to our hypothesis on $y$.
\end{proof}

\begin{proposition}\label{prop:cycle}
Let $p>5$ and $k\leq p-1$.
If $V(k-1)$ exists, then $v_k$ is a permanent cycle in 
the spectral sequence for $V(k-1)$.
\end{proposition}

\begin{proof}
The cases $k=1, 2$ and 3 are immediate, because $v_1$, $v_2$, and $v_3$ are
permanent cycles in the corresponding Adams-Novikov spectral
sequences.  So let $k\geq 4$.  There are many steps to the proof, so
we summarize the argument first.

Consider the connecting homomorphisms
\begin{equation*}
H^*(G;E_{n*}V(i))\overset{\delta_{i+1}}{\lra} H^{*+1}(G;E_{n*}V(i-1))
\end{equation*}
arising from the cofiber sequences 
\begin{equation*}
\Sigma^{|v_{i}|}V(i-1)\overset{v_i}{\lra} V(i-1)\longrightarrow V(i)
\overset{h_{i+1}}{\lra} \Sigma^{|v_i|+1}V(i-1).
\end{equation*}
Let
\begin{equation*}
S^0\overset{f}{\lra} V(1)\overset{g}{\lra} V(k-2)
\end{equation*}
factor the inclusion of the bottom cell into $V(k-2)$.
We will show that 
\begin{equation*}
\delta_k(v_k)=g_*(f_*(\Delta^m)\delta_3(v_3)),
\end{equation*}
where $m\equiv0$ mod~$p$, and that this element is nontrivial.
Since $\Delta^m$ and $v_3$ are permanent
cycles in the spectral sequences for $S^0$ and $V(2)$, respectively,
and $V(1)$ is a ring spectrum~\cite{smi:onreal},
it follows that $\delta_k(v_k)$
is a permanent cycle.

So $\delta_k(v_k)$ detects an element $x\in\pi_{|v_k|-|v_{k-1}|-1}(
E_n^{hG}\wedge V(k-2))$.  By Lemma~\ref{lem:homotopy},
$\pi_{|v_k|-1}(E_n^{hG}\wedge V(k-2))=0$, so $x$ is in the image of
the map
\begin{equation*}
\pi_{|v_k|}(E_n^{hG}\wedge V(k-1))\overset{h_k}{\lra}
\pi_{|v_k|-|v_{k-1}|-1}(E_n^{hG}\wedge V(k-2)),
\end{equation*}
say $x=h_k(y)$.
Note that $y$ must have filtration zero,
i.e., $y\in H^0(G;E_{n*}V(k-1))$.

Now $\delta_k(y-v_k)=0$, so $y-v_k$ pulls back to
$H^0(G;E_{n*}V(k-2))$.  We then show that the map $(g\circ f)_*$ is
onto in this bidegree and that everything in $H^{0,|v_k|}(G;E_{n*})$ is a
permanent cycle.  Therefore $y-v_k$ is a permanent cycle, so $v_k$ is
a permanent cycle.  \smallskip

{\bf Step 1:} $\delta_k(v_k)\ne0$.  Otherwise, $v_k$ is in the image of the
map
\begin{equation*}
H^0(G;E_{n*}V(k-2))\lra H^0(G;E_{n*}V(k-1)),
\end{equation*}
i.e., $v_k$ lifts to a fixed point $\gamma\in u^s\mathbb
F_{q}[\nspc[ u_{k-1},\dots,u_{n-1}]\nspc]$, where $s=1-p^k$.  If
$\mu\in u^s\mathbb F_{q}[\nspc[ u_{k-1},\dots,u_{n-1}]\nspc]$ is a monomial,
we write $\mu\in\gamma$ if $\mu$ appears as a term (up to
multiplication by a unit) when we express $\gamma$ as a sum of
monomials.

We claim that $u^s u_k\not\in\gamma$,
which would then contradict $\delta_k(v_k)=0$.
From~\eqref{eq:action} we have the formula
\begin{equation*}
(\sigma-1)(u^s u_i)\equiv c_i u^su_{i-1}\mod (p,u_1,\dots,u_{i-2},\mathfrak m^2).
\end{equation*}
Suppose $u^s\in\gamma$.  Since $(\sigma-1)\gamma=0$ and $u^su_{n-1}\in
(\sigma-1)u^s$, there must be another monomial $\mu\in\gamma$ with
$u^su_{n-1}\in(\sigma-1)\mu$.  But by the formula above, this isn't possible,
hence $u^s\not\in\gamma$.  Iterating this procedure yields $u^su_i\not\in
\gamma$ for $k\leq i\leq n$.

\smallskip

{\bf Step 2:} The map
\begin{equation*}
H^{1,|v_k|-|v_{k-1}|}(G;E_{n*}V(1))\overset{g_*}{\lra}
 H^{1,|v_k|-|v_{k-1}|}(G;E_{n*}V(k-2))
\end{equation*}
is onto.  It suffices to show that
\begin{equation*}
H^{2,|v_k|-|v_{k-1}|-|v_m|}(G;E_{n*}V(m-1))=0
\end{equation*}
for $m=2,\dots, k-2$, which follows from Lemma~\ref{lem:cohomology}.

\smallskip

{\bf Step 3:} By Step 2, there is an element $y\in H^1(G;E_{n*}V(1))$
which maps to $\delta_k(v_k)\in H^1(G;E_{n*}V(k-2))$.
Then $y=f_*(\Delta^m)\delta_3(v_3)$,
where $m\equiv0$ mod~$p$.  This follows from Proposition~\ref{prop:v1} after 
checking that the degrees work out, i.e.,
\begin{equation*}
|\delta_k(v_k)|\equiv |\delta_3(v_3)|\equiv 2pn+2pn^2\mod 2p^2 n^2,
\end{equation*}
and showing that $\delta_3(v_3)\ne0$, which follows exactly as in Step~1.

\smallskip

 \smallskip
{\bf Step 4:} The map
\begin{equation*}
(g\circ f)_*:H^{0,|v_k|}(G;E_{n*})\longrightarrow H^{0,|v_k|}(G;E_{n*}V(k-2))
\end{equation*}
is surjective.  For this, it suffices to show that
\begin{equation*}
H^{1,|v_k|-|v_m|}(G;E_{n*}V(m-1))=0
\end{equation*}
for $m=0$, 1,\dots,~$k-2$, which follows from Lemma~\ref{lem:cohomology}.

\smallskip

{\bf Step 5:} Everything in $H^{0,|v_k|}(G;E_{n*})$ is a permanent
cycle.  Indeed, for degree reasons these elements are in the kernel $K$
of the projection $H^{0}(G;E_{n*})\to\mathbb F_{q}[\Delta^{\pm}]$.
As we discussed in the previous section, everything in $K$ is a
permanent cycle.
\end{proof}

\begin{proposition}\label{prop:support}
Let $m=(p+3)/2$.  If $V(m-1)$ exists, then $v_m^2$ supports a differential
in the spectral sequence for $V(m-1)$.
\end{proposition}

\begin{proof}
By Lemma~\ref{lem:homotopy},
\begin{equation*}
\pi_{2|v_m|-|v_k|-1}(E_n^{hG}\wedge V(k-1))=0
\end{equation*}
for $k=3$, 4,\dots,~$m-1$.  It follows that if $v_m^2$ is a permanent cycle,
then $v_m^2$ is in the image of the map
\begin{equation*}
H^{0,2|v_m|}(G;E_{n*}V(2))\longrightarrow H^{0,2|v_m|}(G;E_{n*}V(m-1)).
\end{equation*}
We complete the proof by showing that this is not possible.

Set $s=2(1-p^m)$, i.e., $v_m^2=u^s u_m^2$.
Suppose $\gamma\in H^{0,-2s}(G;E_{n*}V(2))$ and let
\begin{equation*}
A_n = \{\,u^s u_3,u^s u_4,\dots,u^s u_{n-1},u^s\,\}
\end{equation*}
and
\begin{equation*}
A_j = \{\,u^s u_i u_j : 2m-j\leq i\leq j\,\}
\end{equation*}
for $j=m$, $m+1$,\dots,~$n-1$.  In the language of Step~1 of the previous
proof, we claim that each of these sets is disjoint from $\gamma$.
This will complete the proof, as
$A_m = \{u^s u_m^2\}$.  This is precisely where the hypothesis
$m=(p+3)/2$ appears.

As in the previous proof, we have
$u^su_i\not\in\gamma$ for $4\leq i\leq n$.  For $i=3$, note that
$\tau(u^su_3)\equiv \eta^{s-1+p^3}u^su_3$ mod $(p,\mathfrak m^2)$,
by~\eqref{eq:action}.  It
is easily verified that $s-1+p^3\equiv(1-p)$ mod $(p-1)^2$, which 
shows that $u^su_3\not\in\gamma$, again by~\eqref{eq:action}.  
(Recall that $\eta$ is an $n^{2}$ root of unity.)
Thus the claim is established for 
$j=n$.

Now let $j<n$.  For each $i$ with $2m-j-1\leq i\leq j$, set
\begin{equation*}
B_{j,i}=\{\,u^su_ku_j\in A_j:k> i\,\}.
\end{equation*}
From~\eqref{eq:action} we have $u^su_{i-1}u_j\in(\sigma-1)u^su_iu_j$.
Then if $u^su_iu_j\in\gamma$, there must be another monomial
$\mu\in\gamma$ with $u^su_{i-1}u_j\in(\sigma-1)\mu$.  We claim that
then 
\begin{equation*}
\mu\in B_{j,i}\cup A_{j+1}\cup\dots\cup A_n.
\end{equation*}
To see this,
write $\mu=u^su_au_b$ with $a\leq b$.  It follows easily
from~\eqref{eq:action} that $b\geq j$ and $a\geq i-1$.  If
$b>j$, then
\begin{equation*}
2m-b\leq 2m-j-1\leq i-1\leq a,
\end{equation*}
hence $\mu\in A_b$.  If $b=j$, then in fact from~\eqref{eq:action}
we have $a\geq i$.  The 
case $a=i$ is ruled out by the assumption that $\mu\ne u^su_iu_j$,
so $\mu\in B_{j,i}$, as claimed.

Suppose inductively that 
$A_k$ is disjoint from $\gamma$ for $k=j+1$,\dots,~$n$.
Suppose also that $B_{j,i}\cap\gamma=\emptyset$ for some $i$
with $2m-j\leq i\leq j$.  (This is trivially true
when $i=j$, as $B_{j,j}=\emptyset$.)  Now by the
previous paragraph, if $u^su_iu_j\in\gamma$, there exists
$\mu\in\gamma\cap(B_{j,i}\cup A_{j+1}\cup\dots\cup A_n)$.  This 
is a contradiction, so $u^su_iu_j\not\in\gamma$, i.e.,
$B_{j,i-1}\cap\gamma=\emptyset$.  By downward induction on $i$,
it follows that $A_j$ is disjoint from~$\gamma$.  By downward 
induction on $j$, we obtain $A_{m}\cap\gamma=\emptyset$, as desired.
\end{proof}

\bibliographystyle{amsalpha}
\bibliography{master}

\end{document}